\documentclass[a4]{article}
\usepackage{amsmath,amsthm,amssymb,graphics}
\usepackage[latin1]{inputenc}
\usepackage{epsfig}
\usepackage[all]{xy}
\usepackage{tikz}

\newtheorem{theorem}{Theorem}
\newtheorem{cor}[theorem]{Corollary}

\begin{document}

\title{New Identities from a Combinatorial Approach to Generalized Fibonacci and Generalized Lucas Numbers}

\author{Robson da Silva\footnote{Supported by CNPq}  \\
	Federal University of S\~ao Paulo - UNIFESP \\
	Department of Science and Technology \\
	S\~ao Jos\'e dos Campos - SP, 12247-014, Brazil \\ silva.robson@unifesp.br \vspace{0.25cm}}

\date{}
\maketitle

\begin{abstract}
We present here some new identities for generalizations of Fibonacci and Lucas numbers by combinatorially interpreting these numbers in terms of numbers of certain tilings of a $1 \times m$ board. As a consequence, some new interesting identities involving the ordinaries Fibonacci and Lucas numbers are derived.
\end{abstract}

\noindent {\bf keyword}: Generalized Fibonacci number, Generalized Lucas number, Tiling


\section{Introduction} 
\label{sect1}

In \cite{Kwasnik} the generalized Fibonacci numbers, $F(k, n)$, and generalized Lucas numbers, $L(k, n)$, were introduced and it was shown that these numbers satisfy, for any integers $k \geq 2$ and $n \geq 0$, the following recurrences:
\begin{equation}
	F(k, n) = \left\{ \begin{array}{l}  n+1, \mbox{ for } n = 0, 1, \ldots, k-1 \\
		F(k, n-1) + F(k, n-k), \mbox{ for } n \geq k
	\end{array} \right.
	\label{Fibo1}
\end{equation}
\begin{equation}
	L(k, n) = \left\{ \begin{array}{l} L(k, n) = n+1, \mbox{ for } n = 0, 1, \ldots, 2k-1 \\
		(k\!-\!1)F(k, n\!-\!(2k\!-\!1))\! +\! F(k, n\!-\!(k\!-\!1)), \mbox{ for } n \geq 2k.
	\end{array}  \right.
	\label{Lucas1}
\end{equation}

For $n \geq 0$, we have that $F(2,n) = F_n$, the $n$th Fibonacci number, and, for $n \geq 3$, $L(2,n) = L_n$, the $n$th Lucas number. Table \ref{tab1} exhibits some values of $F(k, n), L(k, n), F_n$, and $L_n$.


\begin{table}[h]
	\scriptsize
	\centering
	\begin{tabular}{|l|p{0.35cm}|p{0.35cm}|p{0.35cm}|p{0.35cm}|p{0.35cm}|p{0.35cm}|p{0.35cm}|p{0.35cm}|p{0.35cm}|p{0.35cm}|p{0.35cm}|p{0.35cm}|}
		\hline
		$n$ & 0 & 1 & 2 & 3 & 4 & 5 & 6 & 7 & 8 & 9 & 10 & 11 \\ \hline \hline
		$F_n$ & 1 & 2 & 3 & 5 & 8 & 13 & 21 & 34 & 55 & 89 & 144 & 233 \\ \hline
		$F(3,n)$ & 1 & 2 & 3 & 4 & 6 & 9 & 13 & 19 & 28 & 41 & 60 & 88 \\ \hline
		$F(4,n)$ & 1 & 2 & 3 & 4 & 5 & 7 & 10 & 14 & 19 & 26 & 36 & 50 \\ \hline
		$L_n$ & 1 & 2 & 3 & 4 & 7 & 11 & 18 & 29 & 47 & 76 & 123 & 199 \\ \hline
		$L(3,n)$ & 1 & 2 & 3 & 4 & 5 & 6 & 10 & 15 & 21 & 31 & 46 & 67 \\ \hline
		$L(4,n)$ & 1 & 2 & 3 & 4 & 5 & 6 & 7 & 8 & 13 & 19 & 26 & 34 \\ \hline
	\end{tabular}
	\label{tab1}
	\caption{Some values of $F(k, n), L(k, n), F_n$, and $L_n$}
\end{table}

In \cite{Kwasnik} the numbers $F(k, n)$ and $L(k, n)$ are combinatorially interpreted as counting $k$-independent sets of certain finite, undirected, connected, simple graphs, where a subset $A \subset V(G)$ is a $k$-independent set of the graph $G$ if for any two vertices $u, v \in A$, $d_{G}(u,v) \geq k$. For $n \geq 1$, $F(k,n)$ is equal to the number of $k$-independent sets in the graph path or order $n$, $P_n$, and, for $n \geq 3$, $L(k,n)$ is equal to the number of $k$-independent sets in the graph cycle on $n$ vertices, $C_n$. In \cite{Wloch1}, the numbers $F(k,n)$ and $L(k,n)$ are interpreted in terms of the total number of $K_p$-matching in certain graphs.

Here we combinatorially interpret the numbers $F(k, n)$ and $L(k, n)$ in terms of numbers of certain tilings of a $1 \times m$ board, where $m$ is a positive integer. As a consequence some new identities are derived for both the generalized and the ordinaries Fibonacci and Lucas numbers.

Many authors have studied other types of generalizations of the Fibonacci and Lucas numbers, see for instance \cite{Bednarz, Brod, Tasci,Uslu}.

The technique of counting via tilings in different contexts, like in \cite{Benjamin, Benjamin1, Briggs, Little, Little1, Sellers, Stabel}. In \cite{Benjamin1}, for instance, combinatorial interpretations of the Fibonacci and Lucas numbers in terms of certain tilings are used to obtain combinatorial proofs for many identities involving these numbers. Due to the relations \eqref{Fibo1} and \eqref{Lucas1}, the kind of tilings we employ here are different from the ones seen until now.

In \cite{Wloch}, seven identities involving the numbers $F(k,n)$ and $L(k,n)$ are presented. From these identities, when one takes $k=2$, only known identities involving ordinaries Fibonacci and Lucas numbers are found. The main focus of this paper is to exhibit some new identities, 
expanding the list of identities in \cite{Wloch}, which can be useful for counting $k$-independent sets in graphs, and providing, when $k=2$, some new interesting identities involving the usual Fibonacci and Lucas numbers. 


\section{Combinatorial interpretations of $F(k,n)$ and $L(k,n)$}
\label{sect2}

\begin{figure}[h!]
	\centering
	\begin{tikzpicture}[scale=0.6][rounded corners, ultra thick]
	\draw (-0.25,0.25) node {\small $f(3,0) = 1$};
	\shade[top color=black,bottom color=black, draw=black] (1.5,0) rectangle +(0.5,0.5);
	\draw (-0.25,-0.75) node {\small $f(3,1) = 2$};
	\shade[top color=black,bottom color=black, draw=black] (1.5,-1) rectangle +(0.5,0.5);
	\shade[top color=white,bottom color=white, draw=black] (2,-1) rectangle +(0.5,0.5);
	\shade[top color=white,bottom color=white, draw=black] (3,-1) rectangle +(0.5,0.5);
	\shade[top color=black,bottom color=black, draw=black] (3.5,-1) rectangle +(0.5,0.5);
	\draw (-0.25,-1.75) node {\small $f(3,2) = 3$};
	\shade[top color=black,bottom color=black, draw=black] (1.5,-2) rectangle +(0.5,0.5);
	\shade[top color=white,bottom color=white, draw=black] (2,-2) rectangle +(0.5,0.5);
	\shade[top color=white,bottom color=white, draw=black] (2.5,-2) rectangle +(0.5,0.5);
	\shade[top color=white,bottom color=white, draw=black] (3.5,-2) rectangle +(0.5,0.5);
	\shade[top color=black,bottom color=black, draw=black] (4,-2) rectangle +(0.5,0.5);
	\shade[top color=white,bottom color=white, draw=black] (4.5,-2) rectangle +(0.5,0.5);
	\shade[top color=white,bottom color=white, draw=black] (5.5,-2) rectangle +(0.5,0.5);
	\shade[top color=white,bottom color=white, draw=black] (6,-2) rectangle +(0.5,0.5);
	\shade[top color=black,bottom color=black, draw=black] (6.5,-2) rectangle +(0.5,0.5);
	\draw (-0.25,-2.75) node {\small $f(3,3) = 4$};
	\shade[top color=black,bottom color=black, draw=black] (1.5,-3) rectangle +(0.5,0.5);
	\shade[top color=white,bottom color=white, draw=black] (2,-3) rectangle +(0.5,0.5);
	\shade[top color=white,bottom color=white, draw=black] (2.5,-3) rectangle +(0.5,0.5);
	\shade[top color=white,bottom color=white, draw=black] (3,-3) rectangle +(0.5,0.5);
	\shade[top color=white,bottom color=white, draw=black] (4,-3) rectangle +(0.5,0.5);
	\shade[top color=black,bottom color=black, draw=black] (4.5,-3) rectangle +(0.5,0.5);
	\shade[top color=white,bottom color=white, draw=black] (5,-3) rectangle +(0.5,0.5);
	\shade[top color=white,bottom color=white, draw=black] (5.5,-3) rectangle +(0.5,0.5);
	\shade[top color=white,bottom color=white, draw=black] (6.5,-3) rectangle +(0.5,0.5);
	\shade[top color=white,bottom color=white, draw=black] (7,-3) rectangle +(0.5,0.5);
	\shade[top color=black,bottom color=black, draw=black] (7.5,-3) rectangle +(0.5,0.5);
	\shade[top color=white,bottom color=white, draw=black] (8,-3) rectangle +(0.5,0.5);
	\shade[top color=black,bottom color=black, draw=black] (9,-3) rectangle +(0.5,0.5);
	\shade[top color=gray,bottom color=gray, draw=black] (9.5,-3) rectangle +(1.5,0.5);
	\draw (-0.25,-3.75) node {\small $f(3,4) = 6$};
	\shade[top color=black,bottom color=black, draw=black] (1.5,-4) rectangle +(0.5,0.5);
	\shade[top color=white,bottom color=white, draw=black] (2,-4) rectangle +(0.5,0.5);
	\shade[top color=white,bottom color=white, draw=black] (2.5,-4) rectangle +(0.5,0.5);
	\shade[top color=white,bottom color=white, draw=black] (3,-4) rectangle +(0.5,0.5);
	\shade[top color=white,bottom color=white, draw=black] (3.5,-4) rectangle +(0.5,0.5);
	\shade[top color=white,bottom color=white, draw=black] (4.5,-4) rectangle +(0.5,0.5);
	\shade[top color=black,bottom color=black, draw=black] (5,-4) rectangle +(0.5,0.5);
	\shade[top color=white,bottom color=white, draw=black] (5.5,-4) rectangle +(0.5,0.5);
	\shade[top color=white,bottom color=white, draw=black] (6,-4) rectangle +(0.5,0.5);
	\shade[top color=white,bottom color=white, draw=black] (6.5,-4) rectangle +(0.5,0.5);
	\shade[top color=white,bottom color=white, draw=black] (7.5,-4) rectangle +(0.5,0.5);
	\shade[top color=white,bottom color=white, draw=black] (8,-4) rectangle +(0.5,0.5);
	\shade[top color=black,bottom color=black, draw=black] (8.5,-4) rectangle +(0.5,0.5);
	\shade[top color=white,bottom color=white, draw=black] (9,-4) rectangle +(0.5,0.5);
	\shade[top color=white,bottom color=white, draw=black] (9.5,-4) rectangle +(0.5,0.5);
	\shade[top color=black,bottom color=black, draw=black] (1.5,-4.75) rectangle +(0.5,0.5);
	\shade[top color=gray,bottom color=gray, draw=black] (2,-4.75) rectangle +(1.5,0.5);
	\shade[top color=white,bottom color=white, draw=black] (3.5,-4.75) rectangle +(0.5,0.5);
	\shade[top color=black,bottom color=black, draw=black] (4.5,-4.75) rectangle +(0.5,0.5);
	\shade[top color=white,bottom color=white, draw=black] (5,-4.75) rectangle +(0.5,0.5);
	\shade[top color=gray,bottom color=gray, draw=black] (5.5,-4.75) rectangle +(1.5,0.5);
	\shade[top color=white,bottom color=white, draw=black] (7.5,-4.75) rectangle +(0.5,0.5);
	\shade[top color=black,bottom color=black, draw=black] (8,-4.75) rectangle +(0.5,0.5);
	\shade[top color=gray,bottom color=gray, draw=black] (8.5,-4.75) rectangle +(1.5,0.5);
	\draw (-0.25,-5.5) node {\small $f(3,5) = 9$};
	\shade[top color=black,bottom color=black, draw=black] (1.5,-5.75) rectangle +(0.5,0.5);
	\shade[top color=white,bottom color=white, draw=black] (2,-5.75) rectangle +(0.5,0.5);
	\shade[top color=white,bottom color=white, draw=black] (2.5,-5.75) rectangle +(0.5,0.5);
	\shade[top color=white,bottom color=white, draw=black] (3,-5.75) rectangle +(0.5,0.5);
	\shade[top color=white,bottom color=white, draw=black] (3.5,-5.75) rectangle +(0.5,0.5);
	\shade[top color=white,bottom color=white, draw=black] (4,-5.75) rectangle +(0.5,0.5);
	\shade[top color=white,bottom color=white, draw=black] (5,-5.75) rectangle +(0.5,0.5);
	\shade[top color=black,bottom color=black, draw=black] (5.5,-5.75) rectangle +(0.5,0.5);
	\shade[top color=white,bottom color=white, draw=black] (6,-5.75) rectangle +(0.5,0.5);
	\shade[top color=white,bottom color=white, draw=black] (6.5,-5.75) rectangle +(0.5,0.5);
	\shade[top color=white,bottom color=white, draw=black] (7,-5.75) rectangle +(0.5,0.5);
	\shade[top color=white,bottom color=white, draw=black] (7.5,-5.75) rectangle +(0.5,0.5);
	\shade[top color=white,bottom color=white, draw=black] (8.5,-5.75) rectangle +(0.5,0.5);
	\shade[top color=white,bottom color=white, draw=black] (9,-5.75) rectangle +(0.5,0.5);
	\shade[top color=black,bottom color=black, draw=black] (9.5,-5.75) rectangle +(0.5,0.5);
	\shade[top color=white,bottom color=white, draw=black] (10,-5.75) rectangle +(0.5,0.5);
	\shade[top color=white,bottom color=white, draw=black] (10.5,-5.75) rectangle +(0.5,0.5);
	\shade[top color=white,bottom color=white, draw=black] (11,-5.75) rectangle +(0.5,0.5);
	\shade[top color=black,bottom color=black, draw=black] (1.5,-6.5) rectangle +(0.5,0.5);
	\shade[top color=gray,bottom color=gray, draw=black] (2,-6.5) rectangle +(1.5,0.5);
	\shade[top color=white,bottom color=white, draw=black] (3.5,-6.5) rectangle +(0.5,0.5);
	\shade[top color=white,bottom color=white, draw=black] (4,-6.5) rectangle +(0.5,0.5);
	\shade[top color=black,bottom color=black, draw=black] (5,-6.5) rectangle +(0.5,0.5);
	\shade[top color=white,bottom color=white, draw=black] (5.5,-6.5) rectangle +(0.5,0.5);
	\shade[top color=gray,bottom color=gray, draw=black] (6,-6.5) rectangle +(1.5,0.5);
	\shade[top color=white,bottom color=white, draw=black] (7.5,-6.5) rectangle +(0.5,0.5);
	\shade[top color=black,bottom color=black, draw=black] (8.5,-6.5) rectangle +(0.5,0.5);
	\shade[top color=white,bottom color=white, draw=black] (9,-6.5) rectangle +(0.5,0.5);
	\shade[top color=white,bottom color=white, draw=black] (9.5,-6.5) rectangle +(0.5,0.5);
	\shade[top color=gray,bottom color=gray, draw=black] (10,-6.5) rectangle +(1.5,0.5);
	\shade[top color=white,bottom color=white, draw=black] (1.5,-7.25) rectangle +(0.5,0.5);
	\shade[top color=black,bottom color=black, draw=black] (2,-7.25) rectangle +(0.5,0.5);
	\shade[top color=gray,bottom color=gray, draw=black] (2.5,-7.25) rectangle +(1.5,0.5);
	\shade[top color=white,bottom color=white, draw=black] (4,-7.25) rectangle +(0.5,0.5);
	\shade[top color=white,bottom color=white, draw=black] (1.5,-7.25) rectangle +(0.5,0.5);
	\shade[top color=black,bottom color=black, draw=black] (2,-7.25) rectangle +(0.5,0.5);
	\shade[top color=gray,bottom color=gray, draw=black] (2.5,-7.25) rectangle +(1.5,0.5);
	\shade[top color=white,bottom color=white, draw=black] (4,-7.25) rectangle +(0.5,0.5);
	\shade[top color=white,bottom color=white, draw=black] (5,-7.25) rectangle +(0.5,0.5);
	\shade[top color=black,bottom color=black, draw=black] (5.5,-7.25) rectangle +(0.5,0.5);
	\shade[top color=white,bottom color=white, draw=black] (6,-7.25) rectangle +(0.5,0.5);
	\shade[top color=gray,bottom color=gray, draw=black] (6.5,-7.25) rectangle +(1.5,0.5);
	\shade[top color=white,bottom color=white, draw=black] (8.5,-7.25) rectangle +(0.5,0.5);
	\shade[top color=white,bottom color=white, draw=black] (9,-7.25) rectangle +(0.5,0.5);
	\shade[top color=black,bottom color=black, draw=black] (9.5,-7.25) rectangle +(0.5,0.5);
	\shade[top color=gray,bottom color=gray, draw=black] (10,-7.25) rectangle +(1.5,0.5);
	\draw (-0.35,-8) node {\small $f(3,6) = 13$};
	\shade[top color=black,bottom color=black, draw=black] (1.5,-8.25) rectangle +(0.5,0.5);
	\shade[top color=white,bottom color=white, draw=black] (2,-8.25) rectangle +(0.5,0.5);
	\shade[top color=white,bottom color=white, draw=black] (2.5,-8.25) rectangle +(0.5,0.5);
	\shade[top color=white,bottom color=white, draw=black] (3,-8.25) rectangle +(0.5,0.5);
	\shade[top color=white,bottom color=white, draw=black] (3.5,-8.25) rectangle +(0.5,0.5);
	\shade[top color=white,bottom color=white, draw=black] (4,-8.25) rectangle +(0.5,0.5);
	\shade[top color=white,bottom color=white, draw=black] (4.5,-8.25) rectangle +(0.5,0.5);
	\shade[top color=white,bottom color=white, draw=black] (5.5,-8.25) rectangle +(0.5,0.5);
	\shade[top color=black,bottom color=black, draw=black] (6,-8.25) rectangle +(0.5,0.5);
	\shade[top color=white,bottom color=white, draw=black] (6.5,-8.25) rectangle +(0.5,0.5);
	\shade[top color=white,bottom color=white, draw=black] (7,-8.25) rectangle +(0.5,0.5);
	\shade[top color=white,bottom color=white, draw=black] (7.5,-8.25) rectangle +(0.5,0.5);
	\shade[top color=white,bottom color=white, draw=black] (8,-8.25) rectangle +(0.5,0.5);
	\shade[top color=white,bottom color=white, draw=black] (8.5,-8.25) rectangle +(0.5,0.5);
	\shade[top color=white,bottom color=white, draw=black] (1.5,-9) rectangle +(0.5,0.5);
	\shade[top color=white,bottom color=white, draw=black] (2,-9) rectangle +(0.5,0.5);
	\shade[top color=black,bottom color=black, draw=black] (2.5,-9) rectangle +(0.5,0.5);
	\shade[top color=white,bottom color=white, draw=black] (3,-9) rectangle +(0.5,0.5);
	\shade[top color=white,bottom color=white, draw=black] (3.5,-9) rectangle +(0.5,0.5);
	\shade[top color=white,bottom color=white, draw=black] (4,-9) rectangle +(0.5,0.5);
	\shade[top color=white,bottom color=white, draw=black] (4.5,-9) rectangle +(0.5,0.5);
	\shade[top color=black,bottom color=black, draw=black] (5.5,-9) rectangle +(0.5,0.5);
	\shade[top color=gray,bottom color=gray, draw=black] (6,-9) rectangle +(1.5,0.5);
	\shade[top color=white,bottom color=white, draw=black] (7.5,-9) rectangle +(0.5,0.5);
	\shade[top color=white,bottom color=white, draw=black] (8,-9) rectangle +(0.5,0.5);
	\shade[top color=white,bottom color=white, draw=black] (8.5,-9) rectangle +(0.5,0.5);
	\shade[top color=black,bottom color=black, draw=black] (1.5,-9.75) rectangle +(0.5,0.5);
	\shade[top color=white,bottom color=white, draw=black] (2,-9.75) rectangle +(0.5,0.5);
	\shade[top color=gray,bottom color=gray, draw=black] (2.5,-9.75) rectangle +(1.5,0.5);
	\shade[top color=white,bottom color=white, draw=black] (4,-9.75) rectangle +(0.5,0.5);
	\shade[top color=white,bottom color=white, draw=black] (4.5,-9.75) rectangle +(0.5,0.5);
	\shade[top color=black,bottom color=black, draw=black] (5.5,-9.75) rectangle +(0.5,0.5);
	\shade[top color=white,bottom color=white, draw=black] (6,-9.75) rectangle +(0.5,0.5);
	\shade[top color=white,bottom color=white, draw=black] (6.5,-9.75) rectangle +(0.5,0.5);
	\shade[top color=gray,bottom color=gray, draw=black] (7,-9.75) rectangle +(1.5,0.5);
	\shade[top color=white,bottom color=white, draw=black] (8.5,-9.75) rectangle +(0.5,0.5);
	\shade[top color=black,bottom color=black, draw=black] (1.5,-10.5) rectangle +(0.5,0.5);
	\shade[top color=white,bottom color=white, draw=black] (2,-10.5) rectangle +(0.5,0.5);
	\shade[top color=white,bottom color=white, draw=black] (2.5,-10.5) rectangle +(0.5,0.5);
	\shade[top color=white,bottom color=white, draw=black] (3,-10.5) rectangle +(0.5,0.5);
	\shade[top color=gray,bottom color=gray, draw=black] (3.5,-10.5) rectangle +(1.5,0.5);
	\shade[top color=white,bottom color=white, draw=black] (5.5,-10.5) rectangle +(0.5,0.5);
	\shade[top color=black,bottom color=black, draw=black] (6,-10.5) rectangle +(0.5,0.5);
	\shade[top color=gray,bottom color=gray, draw=black] (6.5,-10.5) rectangle +(1.5,0.5);
	\shade[top color=white,bottom color=white, draw=black] (8,-10.5) rectangle +(0.5,0.5);
	\shade[top color=white,bottom color=white, draw=black] (8.5,-10.5) rectangle +(0.5,0.5);
	\shade[top color=white,bottom color=white, draw=black] (1.5,-11.25) rectangle +(0.5,0.5);
	\shade[top color=black,bottom color=black, draw=black] (2,-11.25) rectangle +(0.5,0.5);
	\shade[top color=white,bottom color=white, draw=black] (2.5,-11.25) rectangle +(0.5,0.5);
	\shade[top color=gray,bottom color=gray, draw=black] (3,-11.25) rectangle +(1.5,0.5);
	\shade[top color=white,bottom color=white, draw=black] (4.5,-11.25) rectangle +(0.5,0.5);
	\shade[top color=white,bottom color=white, draw=black] (5.5,-11.25) rectangle +(0.5,0.5);
	\shade[top color=black,bottom color=black, draw=black] (6,-11.25) rectangle +(0.5,0.5);
	\shade[top color=white,bottom color=white, draw=black] (6.5,-11.25) rectangle +(0.5,0.5);
	\shade[top color=white,bottom color=white, draw=black] (7,-11.25) rectangle +(0.5,0.5);
	\shade[top color=gray,bottom color=gray, draw=black] (7.5,-11.25) rectangle +(1.5,0.5);
	\shade[top color=white,bottom color=white, draw=black] (1.5,-12) rectangle +(0.5,0.5);
	\shade[top color=white,bottom color=white, draw=black] (2,-12) rectangle +(0.5,0.5);
	\shade[top color=black,bottom color=black, draw=black] (2.5,-12) rectangle +(0.5,0.5);
	\shade[top color=gray,bottom color=gray, draw=black] (3,-12) rectangle +(1.5,0.5);
	\shade[top color=white,bottom color=white, draw=black] (4.5,-12) rectangle +(0.5,0.5);
	\shade[top color=white,bottom color=white, draw=black] (5.5,-12) rectangle +(0.5,0.5);
	\shade[top color=white,bottom color=white, draw=black] (6,-12) rectangle +(0.5,0.5);
	\shade[top color=black,bottom color=black, draw=black] (6.5,-12) rectangle +(0.5,0.5);
	\shade[top color=white,bottom color=white, draw=black] (7,-12) rectangle +(0.5,0.5);
	\shade[top color=gray,bottom color=gray, draw=black] (7.5,-12) rectangle +(1.5,0.5);
	\shade[top color=black,bottom color=black, draw=black] (1.5,-12.75) rectangle +(0.5,0.5);
	\shade[top color=gray,bottom color=gray, draw=black] (2,-12.75) rectangle +(1.5,0.5);
	\shade[top color=gray,bottom color=gray, draw=black] (3.5,-12.75) rectangle +(1.5,0.5);
	\end{tikzpicture}
	\caption{The tilings for $k=3$ and $n= 0, 1, \ldots, 6$}
	\label{Fig1}
\end{figure}
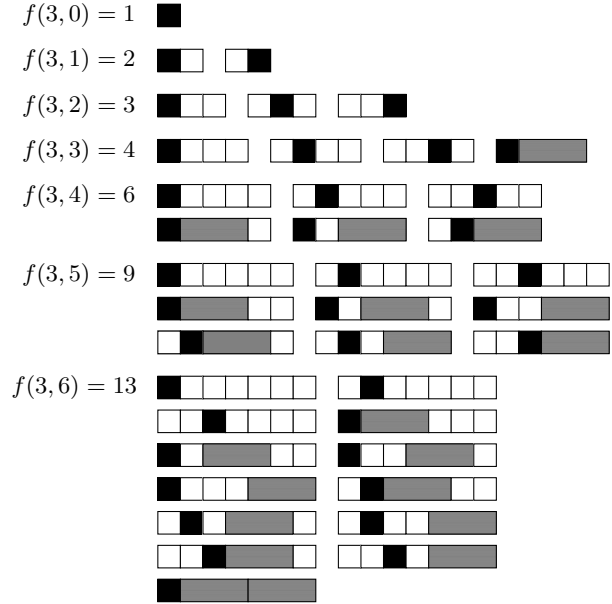

Given integers $k \geq 2$ and $n \geq 0$, we consider tilings of a $1 \times (n+1)$ board using $1 \times 1$ white or black squares and $1 \times k$ gray rectangles such that there is exactly one black square that appears in one of the first $k$ positions. Hence, if $n < k$, there is exactly one black square and all other positions are occupied by white squares. Let $f(k, n)$ be the number of such tilings. Figure \ref{Fig1} shows an example for $k=3$ and $n = 0, 1, \ldots, 6$. 

It is easy to see that for $n < k$, $f(k, n) = n+1$, since there are $n+1$ places to insert the black square, leaving the remaining $n$ spaces filled with white squares. If $n \geq k$, each tiling counted by $f(k,n)$ ends with either a white square or a gray rectangle (never a black square). By removing this last piece we are left with tilings counted either by $f(k, n-1)$ or by $f(k, n-k)$, according to this piece being white or gray, respectively. Then, $f(k,n) = f(k, n-1) + f(k, n-k)$. 

Therefore, as $f(k,n)$ satisfies the same recurrence as $F(k,n)$ and they share the same initial conditions, we have the following combinatorial interpretation for $F(k, n)$:
$$F(k, n) = f(k, n),$$
i.e., $F(k, n)$ is the number of tilings of a $1 \times (n+1)$ board with $1 \times 1$ black or white squares and $1 \times k$ gray rectangles such that there is exactly one black square that appears in one of the first $k$ positions.

As mentioned before, we have $F(2,n) = F_n$, the $n$th Fibonacci number. Hence, when $k=2$, we have a combinatorial interpretation for the Fibonacci numbers. 

From \eqref{Lucas1} and the above interpretation of $F(k,n)$, we can combinatorially interpret $L(k,n)$ as being the number of tilings of a $1 \times (n+1)$ board with $1 \times 1$ white or black squares and $1 \times k$ gray rectangles such that there is exactly one black square that appears in one of the first $k$ positions with an additional  condition:
\begin{itemize}
	\item if $k \leq n < 2k$, there are at least $n-k+1$ pieces after the black square and at least $n-k$ of them are white squares.
	\item if $n \geq 2k$, there are at least $k-1$ white squares among the last $k$ pieces. 
\end{itemize}
In fact, if $n \geq 2k$, those tilings ending with $k-1$ white squares are counted by $F(k, n - (k-1))$ while those whose last $k$ pieces contain exactly one gray rectangle between two of, or after all, the $k-1$ white squares are counted by $(k-1)F(k, n - k - (k-1)) = (k-1)F(k, n - (2k-1))$ since there are $(k-1)$ positions to place the gray piece. When $n < 2k$, it is easy to see that the number of such tilings is $n+1$.

We close this section by pointing out that when $k=2$, the combinatorial interpretation we obtain for the Lucas numbers ($L(2,n) = L_n$) is quite different from that one presented in Chapter 2 of \cite{Benjamin1}.

\section{Main results} 
\label{sect4}

In this section we present new identities that arise from the combinatorial interpretations obtained in Section \ref{sect2}. In what follows, we call type $\mathcal{A}$ and type $\mathcal{B}$ the tilings counted by $F(k,n)$ and $L(k,n)$, respectively.

\begin{theorem}
	Let $k \geq 2$ and $n \geq k$ be integers. Then 
	\begin{equation}
		F(k,n+k) = F(k,n) + \displaystyle\sum_{i=0}^{k-1}F(k, n-i).
		\label{i7}
	\end{equation}
\end{theorem}

\begin{proof}
	The number of tilings enumerated by $F(k,n+k)$ ending with $k$ or more white squares is equal to $F(k, n + k - k) = F(k,n)$. Indeed, removing $k$ of the last white squares we are left with tilings counted by $F(k,n)$. The remaining tilings end with $i$ white squares, where $0 \leq i \leq k-1$. Removing these $i$ white squares together with the last rectangle, we see that the number of such tilings is equal to $F(k, n\!+\!k\! -\! (k\! +\! i)) = F(k, n\! -\! i)$. Hence, adding these numbers we obtain \eqref{i7}.
\end{proof}

Adding $F(k, n-k)$ to both sides of \eqref{i7} we get 
$F(k,n+k) + F(k, n-k) = F(k, n-k) + F(k, n) + \sum_{i=0}^{k-1}F(k, n-i).$
Now, by taking $k=2$, it follows that 
$$F_{n+2} + F_{n-2}  = F_{n-2} + F_n + F_{n} + F_{n-1} = 3F_n,$$ 
which is Identity 7 of \cite{Benjamin1}.

\begin{theorem}
	Let $k \geq 2$ and $n \geq 2k$ be integers. Then 
	\begin{equation}
	F(k,n+1-k) = k +  \displaystyle\sum_{i=0}^{n+1-2k} F(k, n+1 - 2k - i).
	\label{idnova}
	\end{equation}
\end{theorem}

\begin{proof}
Among the tilings counted by $F(k,n+1-k)$, there are $k$ of them with no gray rectangle. We count the remaining tilings according to the position of the rightmost rectangle. By removing this rectangle and the $i$ white squares to its right, we are left with tilings counted by $F(k, n\!+\!1\!-\!k\! -\! k\! -\! i) = F(k, n\!+\!1\!-\!2k\!-\!i)$. Summing over all possible values of $i$ ($0\! \leq\! i\! \leq\! n\!+\!1\!-\!2k$), we obtain \eqref{idnova}.

\end{proof}

Taking $k = 2$, we have the following interesting identity for the Fibonacci numbers.

\begin{cor}
	For $n \geq 4$, we have
	$$\begin{array}{rl}
	F_{n-1} & = 2 + F_{n-3} + F_{n-4} + \cdots + F_{1} + F_{0} \\ 
	& = 2 + \displaystyle\sum_{i=0}^{n-3}F_{i}.
	\end{array}$$
\end{cor}

\begin{theorem}
	Let $k \geq 2$ and $n \geq 2k$ be integers. Then 
	$$ \begin{array}{rl} F(k,n) = & k + k(n+1-k) - \displaystyle\frac{k(k-1)}{2} \\ & +   \displaystyle\sum_{j=0}^{n-2k}\sum_{i=0}^{n-2k-j} F(k, n - 2k - i - j). \end{array}$$
\end{theorem}

\begin{proof}
	It is clear that there are $k$ type $\mathcal{A}$ tilings of a $1 \times (n+1)$ board that do not have gray rectangles. The number of tilings with only one gray rectangle is $k(n+1-k) - \frac{k(k-1)}{2}$. In fact, counting the number of tilings according to the position, $i+1$, of the black square, we have $\sum_{i=0}^{k-1}(n\!+\!1\! -\! k\! -\! i) = k(n\!+\!1\!-\!k)\! -\! \frac{k(k\!-\!1)}{2}$ such tilings.
	
	Now, we count the remaining type $\mathcal{A}$ tilings enumerated by $F(k ,n)$ considering the number of white squares between the last two gray rectangles. There are $F(k, n - 2k - i - j)$ tilings ending with $j$ white squares and having $i$, $i = 0, 1, \ldots, n-2k-j$, white squares between the last two gray rectangles. The maximum value of $j$ is $n-2k$ since there are at least two gray rectangles and, clearly, the values of $i$ depend on the number $j$.
	
	Therefore, summing all these numbers we obtain $F(k, n)$, since we have counted all possible type $\mathcal{A}$ tilings of a $1 \times (n+1)$ board.
\end{proof}

When $k = 2$ we obtain the following identity for the Fibonacci numbers.

\begin{cor}
	For $n \geq 4$, we have
	$$\begin{array}{rl}
	F_n & = 2n-1 + F_{n-4} + 2F_{n-5} + 3F_{n-6} + \cdots \\ & + (n-4)F_{1} + (n-3)F_{0} \\ 
	& = 2n-1 + \displaystyle\sum_{i=0}^{n-4}(i+1)F_{n-4-i}.
	\end{array}$$
\end{cor}

\begin{figure}[h!]
	\center
	\begin{tikzpicture}[scale=0.6][rounded corners, ultra thick]
	\draw (3.4,0.25) node {tail of size $k-1$:};
	\draw (7,0.25) node {$\cdots$};
	\shade[top color=white,bottom color=white, draw=black] (7.5,0) rectangle +(0.5,0.5);
	\draw (8.5,0.25) node {$\cdots$};
	\shade[top color=white,bottom color=white, draw=black] (9,0) rectangle +(0.5,0.5);
	
	\draw (3.875,-0.75) node {tails of size $2k-1$:};
	\draw (7,-0.75) node {$\cdots$};
	\shade[top color=white,bottom color=white, draw=black] (7.5,-1) rectangle +(0.5,0.5);
	\shade[top color=gray,bottom color=gray, draw=black] (8,-1) rectangle +(1.5,0.5);
	\shade[top color=white,bottom color=white, draw=black] (9.5,-1) rectangle +(0.5,0.5);
	\draw (10.5,-0.75) node {$\cdots$};
	\shade[top color=white,bottom color=white, draw=black] (11,-1) rectangle +(0.5,0.5);
	
	\draw (7,-1.75) node {$\cdots$};
	\shade[top color=white,bottom color=white, draw=black] (7.5,-2) rectangle +(0.5,0.5);
	\shade[top color=white,bottom color=white, draw=black] (8,-2) rectangle +(0.5,0.5);
	\shade[top color=gray,bottom color=gray, draw=black] (8.5,-2) rectangle +(1.5,0.5);
	\shade[top color=white,bottom color=white, draw=black] (10,-2) rectangle +(0.5,0.5);
	\draw (11,-1.75) node {$\cdots$};
	\shade[top color=white,bottom color=white, draw=black] (11.5,-2) rectangle +(0.5,0.5);
	
	\draw (7,-2.75) node { $\vdots$};
	
	\draw (7,-3.75) node {$\cdots$};
	\shade[top color=white,bottom color=white, draw=black] (7.5,-4) rectangle +(0.5,0.5);
	\draw (8.5,-3.75) node {$\cdots$};
	\shade[top color=white,bottom color=white, draw=black] (9,-4) rectangle +(0.5,0.5);
	\shade[top color=gray,bottom color=gray, draw=black] (9.5,-4) rectangle +(1.5,0.5);
	\end{tikzpicture}
	\caption{The possible tails for $k=3$}
	\label{Fig2}
\end{figure}
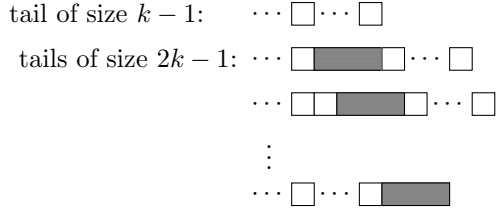

We call \textit{tail} of a type $\mathcal{B}$ tiling the last sequence of $k$ pieces, if a gray rectangle appears, or the last sequence of $k-1$ pieces, if the tiling ends with $k-1$ white squares. Figure \ref{Fig2} shows all possible tails when $k=3$.

The next result provides an interesting recursion for the generalized Lucas numbers and, as particular cases, two recursion for the Lucas numbers, one in terms of Fibonacci numbers and the other involving only Lucas numbers.

\begin{theorem}
	Let $k \geq 2$ and $n \geq 3k-1$ be integers. Then 
	\begin{equation}
		\begin{array}{ll} L(k,n) = & k^2 +  \displaystyle\sum_{i=0}^{n+1-2k}F(k, n - 2k + 1 - i) \\ & +\displaystyle\sum_{i=0}^{n+1-3k}(k-1)F(k, n - 3k + 1 - i). \end{array}
		\label{fim0}
	\end{equation}
\end{theorem}

\begin{proof}
	It is easy to see that there are $k + k(k-1)$ type $\mathcal{B}$ tilings of a $1 \times (n+1)$ board that do not have a gray rectangle before the tail: $k$ having tails of size $k-1$ and $k(k-1)$ having tails of size $2k-1$, since for each position of the black square we have $k-1$ positions to place the gray rectangle in the tail.
	
	In the remaining cases, there is at least one gray rectangle before the tail. We count the number of such tiling by considering the sequence of $i$ white squares between the black square and the leftmost rectangle, where $i \leq n+1 - (k-1) - 1 - k = n-2k+1$ if the tail is of size $k-1$ and $i \leq n+1 - (2k-1) - 1 - k = n-3k+1$ if the tail is of size $2k-1$. Hence, by removing this sequence we are left either with type $\mathcal{A}$ tilings of a $1 \times (n-2k+1-i)$ board, if the tails have size $k-1$ or with type $\mathcal{A}$ tilings of a $1 \times (n-3k+1-i)$ board, if the tails have size $2k-1$. Summing all possible tilings in each case we obtain \ref{fim0}, where the factor $(k-1)$ in the left sum is due to the number of possibilities for placing the gray piece in the tails of size $2k-1$.
\end{proof}

\begin{cor}
	For $n \geq 5$, we have
	\begin{equation}  
		L_n = 4 + F_{n-3} + F_{n-4} + 2\displaystyle\sum_{i=0}^{n-5}F_{n-5-i}
		\label{fim1}
	\end{equation}
	and
	\begin{equation}  
		L_n = 7 + L_{n-2} + L_{n-3} + \cdots + L_{4} + L_{3} = 7 + \displaystyle\sum_{i=0}^{n-5}L_{n-2-i}.
		\label{fim2}
	\end{equation}
\end{cor}

\begin{proof}
	Equation \eqref{fim1} follows immediately from \eqref{fim0} by taking $k=2$. By the recursive definition of $L(k,n)$ given in \ref{Lucas1}, we have $F(k, n-2k+1-i) + (k-1)F(k, n-3k+1-i) \!=\! F(k, (n\!-\!k\!-\!i)\!-\!(k\!-\!1))\! +\! F(k, (n\!-\!k\!-\!i)\!-\!(2k\!-\!1))\! =\! L(k,n\!-\!k\!-\!i)$. Then, 
	$L(k,n)\! =\! k\! +\! k(k\!-\!1)\! +\! F(k,0)\! +\! F(k,1)\! +\! \cdots \!+\! F(k,k\!-\!1)\! +\! \displaystyle\sum_{i=0}^{n-3k+1}L(k,n\!-\!k\!-i\!).$
	Now, taking $k=2$, we obtain \eqref{fim2}.
\end{proof}

Analogous reasoning produces the following identity relating the generalized Lucas numbers and the generalized Fibonacci numbers.

\begin{theorem}
	\small
	Let $k \geq 2$ and $n \geq 4k-1$ be integers. Then 
	$$\begin{array}{rl}
	L(k,n) = & k^2(n+4) - 3k^3 - \frac{k^2(k-1)}{2} \\ & + \! \displaystyle\sum_{j=0}^{n+1-3k}\sum_{i=0}^{n+1-3k-j} F(k, n - 3k + 1 - i - j) \\ & + \! \displaystyle\sum_{j=0}^{n+1-4k}\sum_{i=0}^{n+1-4k-j} (k\!-\!1)F(k, n\! -\! 4k\! +\! 1\! -\! i\! -\! j).
	\end{array}$$
\end{theorem}

\begin{proof}
	We count the number, $L(k,n)$, of type $\mathcal{B}$ tilings of a $1 \times (n+1)$ board, according to the number of gray rectangles before the tails. Clearly there are $k + k(k-1)$ such tilings having no rectangle before the tail: $k$ with tails of size $k-1$ and $k(k-1)$ with tails of size $2k-1$, since for each $k$ possible positions for the black square there are $k-1$ places for the gray rectangle inside the tail.
	
	Now consider those tilings having one gray rectangle before the tail. If the tail is of size $k-1$, for each position $i+1$, $i = 0, 1, \ldots, k-1$, of the black square we have $n+1 - (k-1) - (i+1) - k + 1 = n - 2k + 2 -i$ possible places to insert the gray piece. Hence, the total number of tilings in this case is $\sum_{i=0}^{k-1} (n - 2k + 2 -i) = k(n-2k+2) - \frac{k(k-1)}{2}$. On the other hand, if the tail is of size $2k-1$, for each position $i+1$ of the black square we have $k-1$ for the rectangle in the tail and we are left with $n+1 - (2k-1) - (i+1) - k + 1 = n - 3k + 2 -i$ possible places to insert the gray piece before the tail. The total number of tilings in this last case is $\sum_{i=0}^{k-1} (k-1)(n - 3k + 2 -i) = k(k-1)(n-3k+2) - (k-1)\frac{k(k-1)}{2}$.
	
	Summing all the numbers obtained until now, we have $k^2(n+4) - 3k^3 - \frac{k^2(k-1)}{2}$ type $\mathcal{B}$ tilings of a $1 \times (n+1)$ board having zero or one gray rectangle before the tail.
	
	For the remaining cases, there are at least two gray rectangles before the tails. If we denote by $i$ the number of white squares between the last two gray pieces before the tail and by $j$ the number of white squares between the rightmost rectangle and the tail, we have two cases:
	
	\noindent \textit{Case 1: tail of size $k-1$}. In this case, $i,j \leq n+1 - (k-1) - 1 -2k = n - 3k +1$ and for each $i$ and $j$ fixed the total number of tilings is $F(k, n - 3k + 1 - i - j$. Then, $\sum_{j=0}^{n-3k+1}\sum_{i=0}^{n-3k+1-j}F(k, n-3k+1-i-j)$ counts the number of type $\mathcal{B}$ tilings of a $1 \times (n+1)$ board having tail of size $k-1$ and a sequence of a rectangle, $i$ white squares, a rectangle, and $j$ white squares before the tail. 
	
	\noindent \textit{Case 2: tail of size $2k-1$}. In this case, $i,j \leq n+1 - (2k-1) - 1 -2k = n - 4k +1$ and for each $i$ and $j$ fixed the total number of tilings is $(k-1)F(k, n - 4k + 1 - i - j$, since there are $k-1$ possible places for the gray piece in the tail. Then, $\sum_{j=0}^{n-4k+1}\sum_{i=0}^{n-4k+1-j}(k-1)F(k, n-4k+1-i-j)$ counts the number of type $\mathcal{B}$ tilings of a $1 \times (n+1)$ board having tail of size $2k-1$ and a sequence of a rectangle, $i$ white squares, a rectangle, and $j$ white squares before the tail.
	
\end{proof}

Now, by taking $k=2$, we have the following interesting identity involving the Lucas and Fibonacci numbers.

\begin{cor}
	For $n \geq 7$, we have
	$$\begin{array}{rl}
	L_n & = 4n-10 + F_{n-5} + 2F_{n-6} + 4F_{n-7} + 6F_{n-8} + \\ & + \cdots + (2n-12)F_{1} + (2n-10)F_{0} \\
	& = 4n-10 + F_{n-5} + \displaystyle\sum_{i=1}^{n-5}(2i)F_{n-5-i}
	\end{array}$$
\end{cor}

\section{Concluding remarks}

All identities in the last section can be proved by induction. However the combinatorial proofs explicit the ideas that lead to the identities and, we hope, it can inspire the discovery of new ones.  

The identities involving the numbers $F(k,n)$ and $L(k,n)$, presented in \cite{Wloch}, can be easily proved through the combinatorial interpretations given in Section \ref{sect2}. 

\section{Acknowledgments}
\label{sect5}
This work is supported by the National Council for Scientific and Technological Development (CNPq) of Brazil (No. 473492/2013-0).


\begin{thebibliography}{99}
\bibitem{Benjamin} A. T. Benjamin, S. Plott, J. A. Sellers, Tiling proofs of recent sum identities involving Pell numbers. \textit{Ann. Comb.} \textbf{12} (2008), 271--278.
\bibitem{Benjamin1} A. T. Benjamin, J. J. Quinn, \emph{Proofs That really count: the art of combinatorial proof}. Dolciani Math. Exp., {\bf 27}, 2003.
\bibitem{Bednarz} U. Bednarz, D. Brod, I. W\l och, M. Wolowiec-Musial, On three types of (2,k)-distance Fibonacci numbers and number decompositions. \textit{Ars Comb.} \textbf{118} (2015), 391 -- 405.
\bibitem{Briggs} K. S. Briggs, D. P. Little, J. A. Sellers, Combinatorial proofs of various $q$-Pell identities via tilings. \textit{Ann. Comb.} \textbf{14} (2011), 407--418.
\bibitem{Brod} D. Brod, K. Piejko, I. W\l och, Distance Fibonacci numbers and distance Lucas numbers and their applications. \textit{Ars Comb.} \textbf{112} (2013), 397 -- 410.
\bibitem{Kwasnik} M. Kwasnik, I. W\l och, The total number of generalized stable sets and kernels of graphs. \textit{Ars Combin.} \textbf{55} (2000), 139--146.
\bibitem{Little} D. P. Little, J. A. Sellers, New proofs of identities of Lebesgue and Gollnitz via tilings. \textit{Combin. Theory Ser. A} \textbf{116} (2009), 223--231.
\bibitem{Little1} D. P. Little, J. A. Sellers, A tiling approach to eight identities of Rogers. \textit{European J. of Combin.} \textbf{31} (2010), 694--709.
\bibitem{Sellers} J. A. Sellers, Domino tilings and products of Fibonacci and Pell numbers. \textit{J. Integer Seq.} \textbf{5} (2002), Article 02.1.2.
\bibitem{Stabel} E. C. Stabel, A combinatorial proof of an identity of Ramanujan using tilings. \textit{Bull. Braz. Math. Soc.} \textbf{42} (2011), 203--212.
\bibitem{Tasci} D. Tasci, On the order k generalized Lucas Numbers. \textit{Appl. Math. Comput.} \textbf{155} (2004), 637--641.
\bibitem{Uslu} K. Uslu, N.	Taskara, H. Kose, The Generalized k-Fibonacci and k-Lucas Numbers. \textit{Ars Comb.} \textbf{99} (2011), 25 -- 32.    
\bibitem{Wloch} A. W\l och, Some identities for the generalized Fibonacci numbers and the generalized Lucas numbers. \textit{Appl. Math. Comput.} \textbf{219} (2013), 5564--5568.
\bibitem{Wloch1} A. W\l och, On generalized Fibonacci numbers and $k$-distance $K_p$-matchings in graphs. \textit{Discrete Appl. Math.} \textbf{160} (2012), 1399--1405, 2012.
\end{thebibliography}
\end{document}